\documentclass[10pt,reqno]{amsart}
\usepackage{amsfonts,amsthm,latexsym,amsmath,amssymb,amscd,amsmath,
epsf}
\usepackage{epic}
\usepackage{eepic}
\usepackage{xypic}
\usepackage{graphicx}
\usepackage{mathrsfs}
\xyoption{all}

\theoremstyle{plain}
   \newtheorem{theorem}{Theorem}[section]
   \newtheorem{proposition}[theorem]{Proposition}
   \newtheorem{lemma}[theorem]{Lemma}
   \newtheorem{corollary}[theorem]{Corollary}

     \newtheorem{question}[theorem]{Question}
\theoremstyle{definition}
   
   \newtheorem{example}[theorem]{Example}
\theoremstyle{remark}
\newtheorem{remark}[theorem]{Remark}

\newcommand{\RR}{\mathbb{R}}
\newcommand{\ZZ}{\mathbb{Z}}
\newcommand{\NN}{\mathbb{N}}
\newcommand{\CC}{\mathbb{C}}
\newcommand{\FS}{\mathcal{F}}
\newcommand{\FF}{\mathcal{F}}

\newcommand{\MM}{\mathcal{M}}
\newcommand{\PP}{\mathcal{P}}
\newcommand{\BB}{\mathcal{B}}

\renewcommand{\Im}{\text{Im}}
\renewcommand{\Re}{\text{Re}}
\newcommand{\e}{\ar@{-}}
\newcommand{\de}{\ar@{.}}

\title[Half-plane property and Matroid Theory]{Polynomials with the Half-Plane Property \\ and Matroid theory}

\def\newop#1{\expandafter\def\csname #1\endcsname{\mathop{\rm
#1}\nolimits}}

\newop{St}
\newop{supp}
\newop{diag}
\newop{roots}
\newop{edges}

\author[P.~Br\"and\'en]{Petter Br\"and\'en}
       \address{Department of Mathematics, University of Michigan,
       Ann Arbor, MI 48109-1043, USA}
       \email{branden@umich.edu}

\begin{document}

\maketitle
\nocite{*}\bibliographystyle{plain}
\thispagestyle{empty}

\begin{abstract}
A polynomial $f$ is said to have the half-plane property if there is an open half-plane $H \subset \CC$, whose boundary contains the origin,  
such that $f$ is non-zero whenever all the variables are in $H$. 
This paper answers several open questions relating multivariate polynomials with the half-plane 
property to matroid theory. 
\begin{enumerate}
\item We prove that the support of a multivariate polynomial with the half-plane property is a 
jump system. This answers an open question posed by Choe, Oxley, Sokal and Wagner and generalizes their recent result claiming that the same is true whenever the  polynomial is also homogeneous. 
\item We prove that a multivariate multi-affine polynomial $f \in \RR[z_1,\ldots, z_n]$ has the half-plane property (with respect to the upper half-plane) if and only if 
$$
\frac{\partial f}{\partial z_i}(x)\cdot \frac{\partial f}{\partial z_j}(x) - \frac{\partial^2f}{\partial z_i\partial z_j}(x)\cdot f(x) \geq 0
$$
for all $x \in \RR^n$ and $1 \leq i,j \leq n$. This is used to answer two open questions posed 
by Choe and Wagner regarding strongly Rayleigh matroids. 
\item We prove that the Fano matroid is not the support of a polynomial with the half-plane 
property. This is the first instance of a matroid which does not appear as the support of a polynomial 
with the half-plane property and answers a question posed by Choe et al. 
\end{enumerate}
We also discuss further directions and open problems.
\end{abstract}

\section{Introduction} 
Let $H \subset \CC$ be an open half-plane whose boundary contains the origin. We say that a multivariate polynomial with complex 
coefficients is  $H$-{\em stable} if it is nonzero whenever all the variables are in $H$. Often 
$H= \{ z \in \CC : \Im(z) >0 \}$ or $H=\{z \in \CC : \Re(z) >0 \}$.   If $f$ 
is $H$-stable for some $H$, then $f$ is said to have the {\em half-plane property}. If $H$ is the upper half-plane we say that $f$ is {\em stable}\footnote{There is no standard terminology for the different kinds of stability so our notation differs from some authors.}, and if $H$ is the right half-plane that $f$ is 
{\em Hurwitz stable}. Multivariate polynomials with the half-plane property appear (sometimes hidden) in many different areas such as 
statistical mechanics \cite{heilmann,LY,LS}, complex analysis \cite{hinkkanen,levin}, differential 
equations \cite{BBS2, garding}, engineering \cite{basu, KTM}, optimization \cite{gurvits} and 
combinatorics  \cite{choe,COSW,rayleigh,gurvits,heilmann,wagner1,wagner2}. Recently a striking correspondence between polynomials with the half-plane property 
and matroids was found \cite{COSW}. Choe, Oxley, Sokal and Wagner proved that the 
support of an $H$-stable multi-affine and homogeneous polynomial is the set of bases of a  matroid. A polynomial is multi-affine if it has degree at most one in each variable. The study 
of the relationship between polynomials with the half-plane property and matroid theory 
has since then been continued in a series of papers \cite{choe,rayleigh,gurvits,wagner1,wagner2} where several interesting open questions 
have been raised. In this paper we answer some of these open questions and pose others. 

What if a polynomial with the half-plane property is neither homogeneous, nor multi-affine? What can then 
be said about its support? In \cite{COSW} the problem (Problem 13.3) was raised to find a necessary condition for 
a subset $\FS \subset \NN^n$ to be the support of a polynomial with the half-plane property.   In 
Section \ref{support} we prove that the support of a polynomial with the half-plane property 
is a {\em jump system}. A jump system is a recent generalization of matroids introduced by 
Bouchet and Cunningham \cite{bc1} and further studied by Lov\'asz \cite{lovasz}. This also 
settles Question 13.4 of \cite{COSW}. Prior to this paper no matroids were known {\em not} to be 
the support of a polynomial with the half-plane property and in \cite{COSW} the question (Question 13.7) 
was raised if every matroid is the support of an $H$-stable polynomial. In Section 
\ref{fanosec} we prove 
that the Fano matroid, $F_7$, is not the support of a polynomial with the half-plane property. In Section 
\ref{charac} we prove that a multi-affine polynomial $f \in \RR[z_1, \ldots, z_n]$ is 
stable if and only if 
$$
\frac{\partial f}{\partial z_i}(x)\cdot \frac{\partial f}{\partial z_j}(x) - \frac{\partial^2f}{\partial z_i\partial z_j}(x)\cdot f(x) \geq 0
$$
for all $x \in \RR^n$ and $1 \leq i,j \leq n$. This is used to answer two open questions in \cite{rayleigh}.

\section{Matroids, Delta-Matroids and Jump Systems}
A {\em matroid} is a pair $(\MM, E)$, where $\MM$ is a collection of subsets of a finite set $E$ satisfying, 
\begin{enumerate}
\item $\MM$ is hereditary, i.e., if $B \in \MM$ and $A \subseteq B$, then $A \in \MM$,   
\item The set, $\BB$, of maximal elements with respect to inclusion of $\MM$ respects the {\em exchange axiom}: 
$$
A,B \in \mathcal{B} \mbox{ and } x \in A \setminus B    \Longrightarrow 
$$
$$
\exists y \in B \setminus A \mbox{  such that  } A \setminus \{x\} \cup \{y\} \in \mathcal{B}
$$
\end{enumerate}
The elements of $\MM$ are called {\em independent sets} and the set $\BB$ is called the {\em set of bases of} $\MM$. 
For undefined terminology and more information on matroid theory we refer to \cite{oxley}. 

Bouchet \cite{bouchet1} introduced the notion of a delta-matroid as a generalization of both the independent sets and the set of bases 
of a matroid. A {\em delta-matroid} is a pair $(\FF, E)$, where $\FF$ is a collection of subsets of a finite set $E$ such that 
 $\cup_{A \in \FF}A=E$ and satisfying the following {\em symmetric exchange axiom}: 
 $$
 A,B \in \FF, x \in A \Delta B \Longrightarrow \exists y \in A \Delta B \mbox{ such that } A\Delta \{x,y\} \in \FF.
 $$
 Here $\Delta$ is the {\em symmetric difference} defined by $A \Delta B = (A\cup B) \setminus (A \cap B)$. 
 The independent sets of matroids are precisely those delta-matroids that are hereditary 
 and sets of bases of matroids are precisely the delta-matroids for which all the members of $\FF$ have the same cardinality. 

Jump systems were introduced by Bouchet and Cunningham \cite{bc1} as a generalization of delta-matroids, see also 
\cite{lovasz}. 
Let $\alpha,\beta \in \ZZ^n$ and define $|\alpha|=\sum_{i=1}^n|\alpha_i|$. The set of 
{\em steps from $\alpha$  to $\beta$} is 
defined by 
$$
\St (\alpha,\beta) = \{ \sigma \in \ZZ^n : |\sigma|=1, |\alpha+\sigma-\beta|=|\alpha-\beta|-1\}. 
$$
A collection $\FS$ of points in $\ZZ^n$ is called a {\em jump system} if it respects the following axiom.  

\begin{quote}{\bf Two-step Axiom:} If 
$\alpha,\beta \in \FS$, 
$\sigma \in \St(\alpha,\beta)$ and  $\alpha+\sigma \notin \FS$, then there is 
a $\tau \in \St(\alpha+\sigma,\beta)$ such that $\alpha+\sigma+\tau \in \FS$. 
\end{quote}
Delta-matroids are precisely the jump systems for which $\FF \subseteq \{0,1\}^n$ for some positive integer $n$. For examples of 
matroids, delta-matroids and jump systems see Section \ref{applications}. 

\section{The Support of Polynomials with the Half-Plane Property}\label{support}
An important property of $H$-stable polynomials is that they are closed under taking partial derivatives, see e.g., \cite{BBS2,COSW,LS}. 
\begin{proposition} \label{derivative}
Let $f \in \CC[z_1,\ldots, z_n]$ be $H$-stable. Then either $\partial f / \partial z_1 =0$ or  
$\partial f / \partial z_1 $ is $H$-stable. 
\end{proposition}
If $z_1, \ldots, z_n$ are commuting variables and $\alpha=(\alpha_1,\ldots,\alpha_n) \in \NN^n$ we let 
 $z^\alpha= z_1^{\alpha_1}\cdots z_n^{\alpha_n}$. 
The {\em support}, $\supp(f)$, of a polynomial 
$f(z)= \sum_{\alpha \in \NN^n}a(\alpha)z^\alpha \in \CC[z_1,\ldots, z_n]$ is defined by 
$$
\supp(f) = \{ \alpha \in \NN^n : a(\alpha) \neq 0\}.
$$ 
Equip $\ZZ^n$ with the usual partial order $\leq$, defined by $\alpha \leq \beta$ if $\alpha_i \leq \beta_i$ for all $1 \leq i \leq n$. We write $\alpha < \beta$ if $\alpha \leq \beta$ and $\alpha \neq \beta$. Suppose that $f(z) = \sum_{0\leq \gamma \leq \kappa}a(\gamma)z^\gamma \in \CC[z_1, \ldots, z_n]$ is an $H$-stable polynomial of degree $\kappa_i$ in each variable and suppose that $\alpha, \beta \in \supp(f)$ with $\alpha \leq \beta$. Let 
$$
g(z)= \partial^{\kappa-\beta} [z^\kappa f(1/z)], \ \ \ \ 1/z=(1/z_1,\ldots, 1/z_n). 
$$
It follows that $z^\kappa f(1/z)$ is $\overline{H}$-stable where $\overline{H}=\{\overline{z}: z \in H\}=\{z^{-1}: z \in H\}$, and 
by Proposition \ref{derivative} it follows that  $g(z)$ is also  $\overline{H}$-stable.   
For $\alpha=(\alpha_1, \ldots, \alpha_n) \in \NN^n$ and $f \in \CC[z_1,\ldots, z_n]$ let 
$$
\partial^\alpha f = \frac{ \partial^{\alpha_1}}{\partial z_1^{\alpha_1}} \cdots 
\frac{ \partial^{\alpha_n}}{\partial z_n^{\alpha_n}}f. 
$$
Let 
$$
f_{\alpha,\beta}(z)= \partial^{\alpha}[z^\beta g(1/z)]. 
$$
For $\alpha, \beta \in \ZZ^n$, let  $[\alpha,\beta]= \{\gamma \in \ZZ^n : \alpha \leq \gamma \leq \beta \}$ and  $(\alpha,\beta)= \{\gamma \in \ZZ^n : \alpha < \gamma < \beta \}$. Again, $z^\beta g(1/z)$ is 
$H$-stable, so 
by Proposition \ref{derivative} it follows that $f_{\alpha,\beta}$ is ${H}$-stable and  
$$\supp(f_{\alpha,\beta}) = \{ \gamma - \alpha : \gamma \in \supp(f) \cap [\alpha, \beta]\}.$$

The next theorem says that the support of a 
polynomial with the half-plane property is a jump system. This theorem generalizes the main results of 
\cite{COSW} (Theorem 7.1, Corollary 7.3) and \cite{choe}  (Theorem 2) which say that the same is true 
when in addition the polynomial is homogeneous or 
all terms have degree of the same parity, respectively. 
\begin{theorem}\label{jumpsystem}
Suppose that $f$ has the half-plane property. Then the support of $f$ is a jump system. 
\end{theorem}
\begin{proof}
Since every half-plane can be written as $H=\{ e^{i \theta}z : \Re(z) > 0\}$ for some real $\theta$, it follows that $f$ is $H$-stable 
if and only if $f(e^{-i\theta}z_1, \ldots, e^{-i\theta}z_n)$ is Hurwitz stable. Moreover, 
$\supp(f(z)) = \supp(f(e^{-i\theta}z))$, so we may assume that $f$ is Hurwitz stable.  
Consider $\alpha,\beta \in \supp(f)$.  
Let $\mu(z)$ be the change of variables 
$$ 
z_i \mapsto \begin{cases} 
z_i^{-1} \mbox{ if }  \alpha_i > \beta_i, \\ 
z_i \mbox{ otherwise} 
\end{cases}
$$ 
and let $\gamma \in \NN^n$ be sufficiently large so that  $g(z)=z^{ \gamma}f(\mu(z))$ is a polynomial. 
Clearly $f(z)$ is Hurwitz stable if and only if $g(z)$ is. Moreover $\alpha, \beta \in \supp(f)$ are translated to 
$\alpha', \beta' \in \supp(g)$, where $\alpha' \leq \beta'$. It follows that it is no restriction in assuming that 
$\alpha \leq \beta$, when checking the validity of the two-step axiom.    

Suppose that there is a Hurwitz stable polynomial $f$  and $\alpha, \beta \in \supp(f)$ with 
$\alpha \leq \beta $ for which 
the two-step axiom is violated. Also, let $f$ and $\alpha, \beta$ be minimal with respect to 
$|\alpha -\beta|$. Note that if $f$, $\alpha,\beta \in \supp(f)$ constitutes a counterexample then so does 
$f_{\alpha,\beta}$, $0,\beta-\alpha \in \supp(f_{\alpha,\beta})$. Hence  
we may assume that our minimal counterexample is of the form 
$f(z) = \sum_{\gamma}a(\gamma)z^\gamma \in \CC[z_1,\ldots, z_n]$ with 
$a(0)a(\beta) \neq 0$,  $\beta_i>0$ for all $1 \leq i \leq n$  and  $\supp(f) \subseteq [0, \beta]$.

Let $e_1, \ldots, e_n$ be the standard orthonormal basis of $\RR^n$. By symmetry we may assume that $\sigma=e_1$ in the two-step axiom. Then by the failure of the two-step axiom for this counterexample  we have 
$e_1, 2e_1,e_1+e_2, \ldots, e_1+e_n \notin \supp(f)$, see Fig.~\ref{hasse}.  (If $\beta_1=1$ then $2e_1 \notin [0, \beta]$, and since $\supp(f) \subseteq [0, \beta]$ this gives 
$2e_1 \notin  \supp(f)$.) 
\setlength{\unitlength}{14mm}
\newcommand{\p}{\circle*{0.14}}
\newcommand{\cp}{\circle{0.14}}
\begin{figure}\caption{\label{hasse} Supposed minimal counterexample.}
\begin{center}
\begin{picture}(8.8,4.1)
\put(2.6,.5){
\put(1.5,3){\p}\put(1.7,3){$\beta$}
\put(0.3,1){\cp}\put(-0.15,1){$2e_1$}
\dottedline{0.1}(1,1)(1.1,1.3)
\dottedline{0.1}(2.4,1)(2.26,1.3)
\put(1,1){\cp}\put(1.2,1){$e_1+e_2$}
\put(1.1,1.65){empty}
\put(2.4,1){\cp}\put(2.6,1){$e_1+e_n$}
\put(1,0){\cp}\put(0.7,-0.2){$e_1$}
\put(3,0){\p}\put(2.6,0){$e_2$}\put(3.9,0){\p}\put(4.05,0){$e_n$}
\put(3.24,0){$\ldots$}
\path(2,-0.5)(3,0)
\path(2,-0.5)(3.9,0)
\put(2, -0.5){\p}\put(2.1, -0.7){$0$}
\path(2,-0.5)(1,0)
\put(1.26,0.6){$\cdots$}
\path(1,0)(0.3,1)
\path(1,0)(1,1)
\path(1,0)(2.4,1)
\dottedline{0.1}(3,0)(2.9,0.5)
\dottedline{0.1}(3.9,0)(3.8,0.5)
\dottedline{0.1}(1.5,3)(1.35,2.3)
\dottedline{0.1}(1.5,3)(1.8,2.3)
\dottedline{0.1}(1.5,3)(2.2,2.3)
}
\end{picture}
\end{center}
\end{figure}
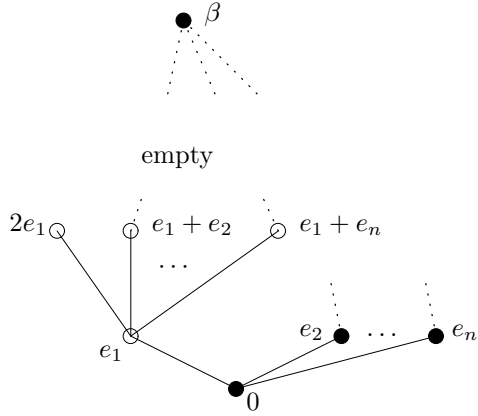
If 
there was a $\xi \in (e_1, \beta)\cap \supp(f)$ then there would be a smaller counterexample $f_{0,\xi}$. Hence, if $\gamma \in \NN^n$ with $\gamma_1>0$ then  $a(\gamma)= 0$ unless $\gamma=\beta$. Let 
$\lambda > 0$ and let $r = 1/{\beta_1}\sum_{i=2}^n\beta_i$. Then the  
univariate polynomial 
$
f(\lambda^{-r}z, \lambda z, \ldots, \lambda z) 
$ 
is Hurwitz stable. Letting $\lambda \rightarrow 0$ we end up with the polynomial  
$$
a(0)+a(\beta)z^{|\beta|},  
$$
which is then Hurwitz stable by Hurwitz's Theorem (on the continuity of the zeros of a polynomial), see 
e.g., \cite[Footnote 3]{COSW} for the appropriate multivariate version.    
We cannot have  $|\beta| \leq 2$, since then the two-step axiom would be valid, so $|\beta| \geq 3$. This gives a contradiction since, when $n \geq 3$, at least one of the 
$n$th roots of a non-zero complex number is in any given half-plane whose boundary contains the origin.    
\end{proof}
An immediate corollary of Theorem \ref{jumpsystem} is a positive answer to Question 13.4 of 
\cite{COSW}. 
\begin{corollary}
The support of a multi-affine polynomial with the half-plane property is a delta-matroid. 
\end{corollary}
 Also Theorem 7.1 of \cite{COSW} follows. 
\begin{corollary}
The support of a multi-affine and homogeneous polynomial with the half-plane property is the 
set of bases of a matroid.
\end{corollary}
\begin{remark}
Recall that the {\em Newton polytope} of a polynomial is the convex hull of its support. In \cite{bc1} it 
was shown that the convex hull of a jump system is a so called {\em bisubmodular polyhedra}, and conversely that the integral points of an integral bisubmodular polyhedra determine a jump system. 
It thus follows that the Newton polytope of a polynomial with the half-plane property is a bisubmodular 
polytope. 
\end{remark}

A polynomial is 
{\em real stable} if it is stable and all coefficients are real. It follows that  a  polynomial 
$f \in \RR[z_1,\ldots, z_n]$ is real stable if and only if  for all lines $z(t) = \lambda t + \alpha$, where 
$\lambda  \in \RR_+^n$ and $\alpha \in  \RR^n$, the polynomial $f(z(t))$ has all zeros 
real. Here $\RR_+$ denotes the set of all positive real numbers. In particular, a univariate polynomial 
with real coefficients is real stable if and only if all its zeros are real. 
\begin{example}
A finite subset $\FF$ of $\NN$ is a jump system if an only if it has holes of size at most $1$, i.e,  
$$
i,k \in \FF, i<k \mbox{ and } j \notin \FF  \mbox{ for all } i<j<k \Longrightarrow k-i \leq 2.
$$ 
Are all finite jump systems in 
$\NN$ supports of polynomials with the half-plane property? Yes! In fact, if we assume that $0 \in \FF$ 
then there is a real-rooted polynomial $f$ with simple zeros such that 
$\FF = \supp(f)$. The proof of this is by induction over the maximal element of $\FF$. If $1 \in \FF$ then 
$$
\FF_1= \{i-1: i \geq 1, i \in \FF \}
$$ 
is a jump system with $0 \in \FF_1$. Hence, by induction, there is a real- and simple-rooted polynomial 
$g$ such that $\supp(g)= \FF_1$. If $\epsilon>0$ is small enough then $
\epsilon + zg$ will be real- and simple-rooted and $\supp(\epsilon +zg) =\FF$. 

If $1 \notin \FF$ then $\FF=\{0\}$ or $2 \in \FF$. In the latter case we have that 
 $$
\FF_2= \{i-2: i \geq 2, i \in \FF \}
$$ 
is a jump system with $0 \in \FF_2$. Hence, by induction,  there is a real- and simple-rooted polynomial 
$g$ such that $\supp(g)= \FF_2$. For small $\epsilon >0$ the polynomial 
$-\epsilon g(0)+z^2g$ will be real- and simple-rooted and $\supp(-\epsilon g(0)+z^2g) =\FF$. 
\end{example}

A well known property 
of real-rooted polynomials with non-negative coefficients is that the coefficients have no internal zeros, 
i.e., if $f(z) = a_0 + a_1z + \cdots + a_nz^n$ is real-rooted and $a_i \geq 0$ for $0 \leq i \leq n$, 
then 
$$
i < j < k \mbox{ and } a_ia_k \neq 0 \Longrightarrow a_j \neq 0.
$$
This extends to several variables:
\begin{corollary} 
Let $f$ be a real stable polynomial with nonnegative coefficients. If $\alpha \leq \gamma \leq \beta$ and  
$\alpha, \beta \in \supp(f)$, then $\gamma \in \supp(f)$. 
\end{corollary}
\begin{proof}
If the corollary is false then there is a real stable polynomial $f$ with nonnegative coefficients,  and points $\alpha, \beta \in \NN^n$ with $\alpha < \beta$, $\alpha, \beta \in \supp(f)$ but  
$\alpha +e_i \notin \supp(f)$ for some $1\leq i \leq n$ with $\alpha +e_i < \beta$. By the two-step axiom there is a $1\leq j \leq n$ such that $\xi = \alpha+e_i+e_j \in \supp(f)$. Now, $f_{\alpha,\xi}=a + bz_j + cz_iz_j$,  $a, b,c \geq 0, ac>0$ is real stable. If $i=j$ then   $f_{\alpha,\xi}=a +  cz_i^2$ is not real stable, so we must have 
$i \neq j$. By letting $z_i =\lambda z$ and $z_j = \lambda^{-1}z$, and letting $\lambda \rightarrow \infty$ 
 we have by Hurwitz's theorem that the univariate polynomial $a+cz^2$ is real stable. This is a contradiction.
\end{proof} 
\section{Applications of the Support Theorem}\label{applications}
Here we give examples of $H$-stable polynomials and their supports. 
\begin{lemma}\label{lincomb} Let $A_i$ be complex positive semidefinite $n \times n$ matrices and let $B$ be complex Hermitian. Then 
$$
f(z)=\det( z_1A_1  + \cdots + z_mA_m + B)
$$
is real stable.
\end{lemma}
\begin{proof}
By Hurwitz's theorem we may assume that the $A_i$'s are all positive definite. Let $z(t) = \lambda t +\alpha$, where $\lambda \in \RR_+^n$ and $\alpha \in \RR^n$. Then $P=\lambda_1A_1+\cdots +\lambda_nA_n$ is positive definite. Thus 
$P$ has a square root, $P^{1/2}$,  and  
$$
f(z(t))= \det(P)\det(tI+P^{-1/2}HP^{-1/2}), 
$$
where $H= B+\alpha_1A_1+\cdots +\alpha_nA_n$ is complex  Hermitian. Hence $f(z(t))$ is a constant multiple of the characteristic polynomial 
of a Hermitian matrix, so all zeros of $f(z(t))$ are real. 
\end{proof} 
In two variables there is a converse to the above lemma, see \cite{BBS2}. 
\begin{theorem}\label{laxlike}Let $f(x,y)\in \RR[x,y]$ be of degree $n$. Then  $f$ is 
real stable if and only if there are two 
$n\times n$ real positive semidefinite matrices $A,B$ and a real symmetric matrix $C$ such that 
$$
f(x,y)=  \pm \det(xA+yB+C).
$$
\end{theorem}
The proof of Theorem \ref{laxlike} uses the Lax Conjecture on hyperbolic polynomials which 
was proved only very recently \cite{LPR}, see also \cite{helton,vinnikov}. 
Let $Z=\diag(z_1, \ldots, z_n)$ be a diagonal matrix. Consequences of Lemma~\ref{lincomb} are the following. 
\begin{itemize} 
\item If $A$ is a Hermitian $n \times n$ matrix then the polynomials $\det(Z+A)$ and $\det(I+AZ)$ are real stable, 
\item If $A$ is a skew-Hermitian $n \times n$ matrix then $\det(Z+A)$ and  $\det(I+AZ)$ are Hurwitz stable. 
\end{itemize}
For an $n \times n$ matrix $A$ let  $A[S]$ denote the  principal  sub-matrix of $A$ with rows and columns in $A$ indexed by $S \subseteq \{1,\ldots, n\}$.  
In \cite{bouchet1} Bouchet proved that  the set 
$$
\{ S  \subseteq \{1,\ldots, n\} : A[S] \mbox{ is non-singular}\}
$$
is a delta-matroid whenever $A$ is a $n \times n$ symmetric or skew-symmetric 
matrix over a field. The proof is not trivial. However, when the field is $\CC$ it follows as a corollary of Theorem \ref{jumpsystem}. 
 
\begin{corollary} 
Let $A$ be a Hermitian or a skew-Hermitian $n \times n$ matrix and let 
$\mathcal{F}= \{ S \subseteq \{1,\ldots, n\} : A[S] \mbox{ is non-singular}\}$, where $A[S]$ is the principal minor 
with rows and columns indexed by $S$. Then $\mathcal{F}$ is a $\Delta$-matroid. 
\end{corollary}
\begin{proof}
The corollary follows from Theorem \ref{jumpsystem} and the fact that 
$$
\det(I+AZ) = \sum_{S \subseteq \{1, \ldots, n\}} \det(A[S])z^S, 
$$
is stable (Hurwitz stable), so $\supp(\det(I+ZA)) = \FF$.
\end{proof}
The general form of the Heilmann-Lieb Theorem \cite{heilmann} is the following. 
\begin{theorem}[Heilmann-Lieb] Let $G=(V,E)$ be a graph,  
$V=\{1,\ldots,n\}$. To each edge $e =ij \in E$ assign a non-negative real number $\lambda_{ij}$.Then 
the polynomial 
$$
M_G(z)=\sum_ {M \mbox{\tiny is a matching }}\prod_{ij \in M}\lambda_{ij}z_iz_j
$$
is Hurwitz-stable. 
\end{theorem}
As a corollary of the Heilmann-Lieb Theorem and Theorem \ref{jumpsystem} we get the 
following result which is usually proved using augmented path arguments. 
\begin{corollary}
Let $G=(V,E)$ be a graph and let $\FF$ be the collection of subsets of $V$ consisting of all $S$ for which there is a matching of $G$ covering 
precisely the elements of $S$. Then $\FF$ is delta-matroid. 
\end{corollary}
\begin{proof}
The corollary follows from the  Heilmann-Lieb Theorem (letting $\lambda_{ij}=1$) and Theorem \ref{jumpsystem} since 
$\supp(M_G(z))=\FF$.
\end{proof}
Let $G=(E,V)$ where $V=\{1,\ldots, n\}$ and let $D(G)=(d_1, \ldots, d_n)$ be the 
{\em degree sequence} of $G$. Here $d_i$ is the degree of the vertex $i$. The 
polynomial 
$$
\sum_{H=(F,V), F \subseteq E} z^{D(H)} = \prod_{ij \in E}(1+z_iz_j)
$$
is clearly Hurwitz stable.  As a consequence of Theorem \ref{jumpsystem} we get that  
$$
\{ D(H) : H \mbox{ spanning subgraph of } G \} 
$$
is a jump system. 

Let $A$ be an $r \times n$ matrix with complex entries and let $A^*$ be its complex adjoint. 
By the Cauchy-Binet formula we have 
 \begin{eqnarray*} 
\det(AZA^{*}) &=& \det\left( \sum_{k=1}^r z_k (a_{ik}\overline{a_{jk}})_{1 \leq i,j \leq r} \right) \\
                        &=& \sum_{S \in \binom {[n]} r } |\det A[S]|^2 z^S. 
\end{eqnarray*}
Since  $(a_{ik}\overline{a_{jk}})_{1 \leq i,j \leq r}$ is positive semi definite, $1\leq k \leq r$,   
Lemma \ref{lincomb} gives that any matroid representable over $\CC$ is the support 
of a real stable polynomial, see \cite{COSW} for another proof.  

Let $G=(V,E)$ be a graph with vertex set $V=\{1, \ldots, n\}$ and edge set $E$. Associate to each edge $e \in E$ a variable $w_e$. If $e$ connects $i$ and $j$ let $A_e$ be the $n \times n$ positive semidefinite matrix with the $ii$-entry and 
the $jj$-entry equal to $1$, with the $ij$-entry and $ji$-entry equal to $-1$ and with all other entries equal to $0$. The {\em Laplacian}, $L(G)$,  of $G$ may be defined by 
$$
L(G)= \sum_{e \in E}w_eA_e.
$$
Let
$$
f_G(z,w)= \det(L(G)+Z).
$$
Thus, by Lemma~\ref{lincomb}, $f_G$ is a multi-affine real stable polynomial with non-negative coefficients.   
The Principal Minors Matrix-Tree Theorem (see e.g. \cite{Chaiken}) says that 
$$
 f_G(z,w)= \sum_{F}z^{\roots(F)}w^{\edges(F)}, 
$$ 
where the sum is over all rooted spanning forests $F$ in $G$, $\roots(F) \subseteq V$ is the set of roots 
of $F$ and $\edges(F) \subseteq E$ is the set of edges used in $F$.  Since the class of stable polynomials is closed under differentiation and specialization of  variables at real values (see e.g. \cite{BBS2}) we have that the {\em spanning tree polynomial} 
$$
T_G(w)= \sum_Tw^T=\frac {\partial f}{\partial z_i}\Big|_{z=0}
$$ 
where the sum is over all spanning trees in $G$ is real stable (which is widely known). The support 
of $T_G(w)$ is the graphic matroid associated with $G$. 

\section{A Characterization of Real Stable Multi-affine Polynomials}\label{charac}
Here we will give a characterization of real stable multi-affine polynomials. First we will need 
some  results on univariate stable polynomials and some results from \cite{BBS2}.   
Let $\alpha_1 \leq \alpha_2 \leq \cdots \leq \alpha_n$ and 
$\beta_1 \leq \beta_2 \leq \cdots \leq \beta_m$ be the zeros of two univariate polynomials with real zeros only. The zeros are 
{\em interlaced} if they can 
be ordered so that $\alpha_1 \leq \beta_1 \leq \alpha_2 \leq \beta_2 \leq 
\cdots$ or $\beta_1 \leq \alpha_1 \leq \beta_2 \leq \alpha_2 \leq 
\cdots$. Note that by our convention, the zeros of any two polynomials 
of degree $0$ or $1$ interlace. It is not hard to see that if 
the zeros of $h$ and $g$ interlace then the {\em Wronskian}, 
$W[g,h]=g'h-gh'$ is either non-negative or non-positive on the whole of 
$\RR$. Let $g,h \in \RR[z]$. We say that $g$ and $h$ are  in  
{\em proper position}, denoted $g \ll h$ if 
the zeros of $h$ and $g$ interlace and  $W[g,h] \leq 0$. 
For technical reasons we also say that the zeros of the polynomial $0$ 
interlaces the zeros of any (non-zero) real-rooted polynomial $f$, and write 
$0 \ll f$ and $f \ll 0$. The Hermite-Biehler Theorem characterizes univariate stable 
polynomials, see \cite{RS}.
\begin{theorem}[Hermite-Biehler]
Let $f = h +ig \in \CC[z]$ where $h,g \in \RR[z]$. Then 
$f$ is stable if and only if $g \ll h$. 
\end{theorem}
Obreschkoff's Theorem describes linear pencils of polynomials with real zeros only, see 
\cite{obreschkoff,RS}. 
\begin{theorem}[Obreschkoff]
Let $g,h \in \RR[z]$. Then all non-zero polynomials in the pencil 
$$
\{ \alpha h + \beta g : \alpha, \beta \in \RR \}
$$
are real-rooted if and only if $h \ll g$, $g \ll h$ or  $h=g=0$. 
\end{theorem}
We extend the notion of proper position to multivariate polynomials as follows. Two multivariate polynomials $g,h \in \RR[z_1, \ldots, z_n]$ are said to be in {\em proper position}, denoted 
$g \ll h$,  if 
\begin{equation}\label{sevdef}
g(\alpha + vt) \ll h(\alpha +vt)
\end{equation}
for all $\alpha \in \RR^n$ and $v \in \RR_+^n$. Note that for univariate polynomials the two definitions of proper position coincide. The Hermite-Biehler Theorem and Obreschkoff's 
Theorem have the following extensions to several variables, see \cite{BBS2}. 
\begin{theorem}\label{HBseveral}
Let $f = h +ig \in \CC[z_1,\ldots, z_n]$ where $h,g \in \RR[z_1,\ldots,z_n]$. Then 
$f$ is stable if and only if $g \ll h$. 
\end{theorem}
\begin{theorem}\label{obreschkoffseveral}
Let $h,g \in \RR[z_1,\ldots,z_n]$. Then all non-zero polynomials in the pencil 
$$
\{ \alpha h + \beta g : \alpha, \beta \in \RR \}
$$
are real stable if and only if $h \ll g$, $g \ll h$ or  $h=g=0$. 
\end{theorem}
By combining the previous two theorems we get.
\begin{corollary}\label{addz}
Let $f=h+ig\neq 0$ where $h,g \in \RR[z_1,\ldots,z_n]$, and let $z_{n+1}$ be a new indeterminate. Then the following are equivalent. 
\begin{itemize}
\item[(a)] $f=h+ig$ is stable, 
\item[(b)] $h+z_{n+1}g$ is real stable, 
\item[(c)] all nonzero polynomials in the pencil 
$$
\{ \alpha h + \beta g : \alpha, \beta \in \RR \}
$$
are real stable and  
$$
\frac {\partial h} {\partial z_j}(x)\cdot g(x)-h(x)\cdot\frac {\partial g}{\partial z_j}(x) \geq 0,  
$$ 
for all $1\leq j \leq n$ and $x \in \RR^n$. 
\end{itemize}
\end{corollary}
\begin{proof}
$(b) \Rightarrow (a)$: If $h+z_{n+1}g$ is real stable, then in particular it is stable. Hence, since $\Im(i) >0$, we have that 
$h+ig$ is stable. 

$(a) \Rightarrow (c)$: If (a) is true then the statement about the pencil in 
(c) follows immediately from 
Theorem \ref{HBseveral} and Theorem \ref{obreschkoffseveral}. Let $v$ be a vector in $\RR_+^n$. Then, by \eqref{sevdef}, 
$
g(x+(e_j+\epsilon v)t) \ll h(x+(e_j+\epsilon v)t)
$
 so $W[g(x+(e_j+\epsilon v)t),  h(x+(e_j+\epsilon v)t)] \leq 0$ for all $1 \leq j \leq n$ and $x \in \RR^n$. Letting $\epsilon \rightarrow 0$ we have by continuity that 
$$
\frac {\partial h} {\partial z_j}(x)\cdot g(x)-h(x)\cdot\frac {\partial g}{\partial z_j}(x) = 
-W(g(x+e_jt), h(x+e_jt))\big|_{t=0} \geq 0. 
$$ 

$(c) \Rightarrow (b)$: Fixing $z_{n+1}= a+ib$, we have to prove that $h+(a+ib)g=(h+ag)+ibg$ is stable whenever $a \in \RR$ and 
$b \in \RR_+$. If $\alpha, \beta \in \RR$ then $\alpha (h+ag)+\beta b g=\alpha h + (a\alpha +b \beta)g$ 
is either real-stable or identically zero by assumption. Since we do not have $bg = h+ag=0$ we have by Theorem \ref{obreschkoffseveral} that 
$h+ag \ll bg$ or $bg \ll h + ag$. Now, 
\begin{eqnarray*}
W(bg(\alpha+vt),h(\alpha+vt)+ag(\alpha+vt)) &=& bW(g(\alpha+vt),h(\alpha+vt)) \\
&=& -b\sum_{j=1}^nv_j\big(\frac {\partial h} {\partial z_j}\cdot g-h\cdot\frac {\partial g}{\partial z_j}\big)(\alpha+vt) \leq 0
\end{eqnarray*}
whenever $\alpha \in \RR^n$, $v \in \RR_+^n$ and $t \in \RR$. The conclusion now follows from Theorem \ref{HBseveral} and \eqref{sevdef}.

\end{proof}
Using this corollary we may characterize real stable multi-affine polynomials as follows. 
For $f \in \CC[z_1, \ldots, z_n]$ and $1 \leq i,j \leq n$ let 
$$\Delta_{ij}(f)= \frac{\partial f}{\partial z_i}\cdot \frac{\partial f}{\partial z_j} - \frac{\partial^2f}{\partial z_i\partial z_j}\cdot f=- f^2 \frac{\partial^2}{\partial z_i \partial z_j} \big[\log|f| \big].
 $$
\begin{theorem}\label{logf}
Let $f \in \RR[z_1,\ldots,z_n]$ be multi-affine. Then the following are equivalent
\begin{itemize}
\item[(1)] For all $x \in \RR^n$ and $1 \leq i,j \leq n$ 
$$
\Delta_{ij}(f)(x) \geq 0,
$$
\item[(2)] $f$ is stable.
\end{itemize}
\end{theorem}
 \begin{proof}
 (2) $\Rightarrow$ (1): Write $f$ as $f=h+z_i g$, where 
 $h=  f|_{z_i=0}$ and $g= \frac{\partial f}{\partial z_i}$. Then 
 \begin{eqnarray*}
 \Delta_{ij}(f)&=& g\cdot (\frac {\partial h}{\partial z_j}+z_i\frac {\partial g}{\partial z_j})-\frac {\partial g}{\partial z_j}\cdot (h+z_i g)\\
 &=& g\cdot \frac {\partial h}{\partial z_j}-\frac {\partial g}{\partial z_j}\cdot h,
 \end{eqnarray*}
 so by Corollary~\ref{addz}, (b) $\Rightarrow$ (c)  this direction follows.

(1) $\Rightarrow$ (2): The proof is by induction over $n$. Write $f \in \RR[z_1, \ldots, z_{n+1}]$ as 
 $f = h+z_{n+1}g$. We want to apply Corollary~\ref{addz}, (c) $\Rightarrow$ (b).  Let $\alpha \in \RR$, $x \in \RR^n$ and   $1 \leq i,j \leq n$. Then 
 $$
 \Delta_{ij}(f\big|_{z_{n+1}=\alpha})(x) = \Delta_{ij}(f)(x_1, \ldots, x_n, \alpha)\geq 0. 
 $$
 By induction $h + \alpha g$ is real stable or $h + \alpha g=0$. This verifies the condition about the pencil in Corollary~\ref{addz} (c).  Also, 
 $$
\frac {\partial h} {\partial z_j}(x)\cdot g(x)-h(x)\cdot\frac {\partial g}{\partial z_j}(x) = 
\Delta_{j,{n+1}}(f)(x)\geq 0
$$ 
for all $x \in \RR^n$ and $1 \leq j \leq n$. This verifies the Wronskian condition in Corollary~\ref{addz} (c)  which completes the proof.
 \end{proof}
The above theorem is not true without the requirement that $f$ is multi-affine. However,  for 
non-multi-affine polynomials it is still true that (2) $\Rightarrow$ (1). 
 \begin{example}\label{n=2}
 Consider $f(z_1,z_2)= a_{00}+ a_{01}z_2+a_{10}z_1+a_{11}z_1z_2 \in \RR[z_1,z_2]$. Then 
 $$
 \Delta_{12}(f)= -
\left| \begin{array}{cc}
 a_{00} & a_{01} \\
a_{10} & a_{11} 
\end{array} \right|, 
$$
so $f$ is real stable if and only if $\det (a_{ij}) \leq 0$. 
\end{example}
\subsection{The non-multi-affine case}
For the non-multi-affine case we may apply the Grace-Walsh-Szeg\"o Coincidence Theorem 
\cite{grace,szego,walsh}.  
Let $f \in \CC[z_1,\ldots, z_n]$ be a polynomial of degree $d_i$ in the variable $z_i$ for $1\leq i \leq n$. The 
{\em polarization}, $\PP(f)$, is the unique polynomial in the variables 
$\{ z_{ij} : 1 \leq i \leq n, 1 \leq j \leq d_i\}$ satisfying 
\begin{enumerate}
\item $\PP(f)$ is multi-affine, 
\item $\PP(f)$ is symmetric in the variables $z_{i1}, \ldots z_{id_i}$ for $1 \leq i \leq n$, 
\item If we let  $z_{ij}=z_i$ for all $i,j$ in $\PP(f)$ we recover $f$.
\end{enumerate}
A {\em circular region} in $\CC$ is either an open or closed affine half-plane or the open or closed interior or exterior of a circle.  
\begin{theorem}[Grace-Walsh-Szeg\"o]
Let $f \in \CC[z_1,\ldots, z_n]$ be symmetric and multi-affine and let $C$ be a circular region 
containing the points $\zeta_1, \ldots, \zeta_n$. Then there exists a point 
$\zeta \in C$ such that 
$$
f(\zeta_1, \ldots, \zeta_n)=f(\zeta, \ldots, \zeta).
$$
\end{theorem}
From the Grace-Walsh-Szeg\"o Theorem we immediately deduce:
\begin{corollary}\label{polarizeit}
Let $f \in \CC[z_1,\ldots, z_n]$ and let $H$ be a half-plane in $\CC$. Then 
$f$ is $H$-stable if and only if $\PP(f)$ is $H$-stable. 
\end{corollary}

\begin{theorem}\label{polchar}
Let $f \in \RR[z_1,\ldots,z_n]$ be of degree $d$. Then the following are equivalent
\begin{enumerate}
\item For all $x \in \RR^d$ and $1 \leq i,j \leq d$ 
$$
\Delta_{ij}(\PP(f))(x) \geq 0,
$$
\item $f$ is stable.
\end{enumerate}
\end{theorem}
\begin{remark} In the univariate case Theorem~\ref{polchar} gives a characterization of 
polynomials with real zeros only. Since the polarization of a univariate polynomial is 
symmetric we get a single equation in $n-2$ variables, where $n$ is the degree of the 
polynomial. This raises the problem of testing polynomial inequalities which are symmetric all variables and of degree at most two in each variable. It would also be interesting to compare this characterization with the classical \cite[p. 203]{gantmacher}. 
\end{remark}
\subsection{Balanced, Rayleigh and HPP Matroids}
Feder and Mihail \cite{feder} introduced the concept of a balanced matroid in relation to a conjecture 
of Mihail and Vazirani \cite{mihail} regarding expansion properties of one-skeletons 
of $\{0,1\}$-polytopes. Let $\mathcal{M}$ be a matroid on a ground-set $E$. For disjoint subsets $I,J$ of 
$E$ let $\mathcal{M}_I^J$ be the minor of $ \mathcal{M}$ obtained by contracting $I$ and deleting 
$J$, 
$$
\mathcal{M}_I^J= \{ S \setminus I  : S \in \mathcal{M},  I \subseteq S \subseteq E\setminus J \}.
$$      
Let $M_I^J$ denote the number of bases of $\mathcal{M}_I^J$ and let 
$$
\Delta_{ij}(\mathcal{M}):=M_i^jM_j^i-M_{ij}M^{ij}.  
$$ 
Feder and Mihail say that 
$\mathcal{M}$ is {\em balanced} if $\Delta_{ij}(\mathcal{M}_I^J)\geq 0$ 
for all $i,j \in E$ and all $I,J \subseteq E$. Motivated by a property of linear resistive electrical networks, Choe and Wagner 
introduced the notion of a {\em Rayleigh matroid}. Let 
$
M(z):= \sum_{B}z^B 
$ 
where the sum is over the set of bases of $\mathcal{M}$. Then $\mathcal{M}$ is said to be 
Rayleigh if 
$$
\Delta_{ij}(M)(z):=M_i^j(z)M_j^i(z)-M_{ij}(z)M^{ij}(z) \geq 0,  
$$ 
for all $z \in \RR_+^n$ and $i,j \in E$ ($|E|=n$).  It follows that  a matroid is balanced if it is Rayleigh. 
A matroid, $\mathcal{M}$, is {\em strongly Rayleigh} if 
$$
\Delta_{ij}(M)(z):=M_i^j(z)M_j^i(z)-M_{ij}(z)M^{ij}(z) \geq 0,  
$$ 
for all $z \in \RR^n$ and $i,j \in E$.  A matroid $\mathcal{M}$ is HPP (half-plane property) if 
$M(z)$ has the half-plane property. It was proved in \cite{rayleigh} that a matroid is Rayleigh 
if it is HPP. We may now answer the following two open questions posed in \cite{rayleigh}. 
\begin{question}[Choe-Wagner]
Are all HPP-matroids strongly Rayleigh? 
\end{question}
\begin{question}[Choe-Wagner]
Are all strongly Rayleigh matroids HPP? 
\end{question}
By Theorem \ref{logf}  and the identity 
$$
\Delta_{ij}(M)(z) = \Delta_{ij}(M(z)) 
$$ 
both of these question have positive answers. 
\begin{corollary} 
A  matroid is strongly Rayleigh if and only if it is HPP.
\end{corollary}

\section{The Fano Matroid}\label{fanosec}
Consider Figure~\ref{fano}. The set of bases of the {\em Fano matroid}, $F_7$, is the 
collection of subsets of $\{1,\ldots, 7\}$ of cardinality $3$ that are not on a line, i.e., 
all subsets of cardinality $3$ except $\{1,2,3\}$, $\{3,4,5\}$, $\{1,5,6\}$, $\{1,4,7\}$, $\{2,5,7\}$, 
$\{3,6,7\}$ and $\{2,4,6\}$. 
 
\setlength{\unitlength}{20mm}
\newcommand{\p}{\circle*{0.11}}
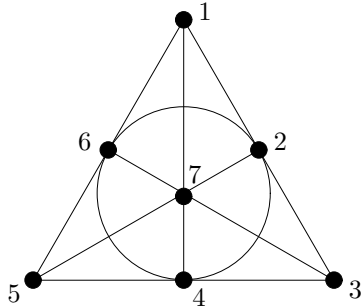
\begin{figure}\caption{\label{fano} The Fano matroid $F_7$.}
\begin{center}
\begin{picture}(17,3)
\put(1.3,0.5){
\put(2,2){\p} \put(2.1,2){1} 
\put(2,0.26795){\p}\put(2.06,0.1){4}
\put(1,0.26795){\p}\put(0.83,0.12){5}
\put(3,0.26795){\p}\put(3.1,0.15){3}
\path(2,2)(2,0.26795)
\path(1,0.26795)(3,0.26795)
\path(2,2)(3,0.26795)
\put(2.5,1.13397){\p}\put(2.6,1.13397){2}
\put(1.5,1.13397){\p}\put(1.3,1.13397){6}
\path(2,2)(1,0.26795)
\path(1,0.26795)(2.5,1.13397)
\path(3,0.26795)(1.5,1.13397)
\put(2, 0.825299){\p}\put(2.03, 0.915){7}
\put(2, 0.846){\circle{1.15}}     
}
\end{picture}
\end{center}
\end{figure}

The Fano matroid is represented over a field of cardinality $2$ by the 
matrix 

$$
 \begin{array}{ccccccc}
1 & 2 & 3 & 4& 5 &6 & 7
\end{array} 
$$
$$
\left[ \begin{array}{ccccccc}
1 & 1 & 0 & 0 & 0 & 1& 1\\
0 & 0 & 0 & 1& 1 & 1& 1 \\
0 & 1 & 1 & 1& 0 &0 & 1  
\end{array} \right]
$$

This configuration of $7$ lines is known as the {\em Fano projective plane}. The Fano matroid has more 
symmetry than Figure~\ref{fano} suggests. The automorphism group of the Fano matroid acts transitively on its point set, and on its line set. (In fact, it acts transitively on the set of ordered triples of non-collinear 
points.)  

We will in this section prove that the Fano matroid is not the support of a polynomial with the 
half-plane property. This is the first instance in the literature of a matroid which is not the support 
of an $H$-stable polynomial and answers Question 13.7 of \cite{COSW}. 
\begin{lemma}\label{key}
Let $f=\sum_{T \subseteq \{1,\ldots,n\}}a(T)z^T$ be a homogeneous multi-affine polynomial with  the half-plane property. Suppose that 
$S\cup \{i,j\} \notin \supp(f)$. Then 
$$
a(S\cup \{i,k\})a(S\cup \{j,\ell\})=a(S\cup \{i,\ell\})a(S\cup \{j,k\})
$$
for all $k, \ell \in \{1,\ldots,n\}$.
\end{lemma}
\begin{proof}
The coefficients of a homogeneous polynomial with the half-plane property all have the same phase, i.e., the 
quotient of two non-zero coefficients is always a positive number, see \cite[Theorem 6.1]{COSW}. Hence 
we may assume that all coefficients are real and non-negative. 

By considering $f_{S, S\cup\{i,j,k,\ell\}}$ we may assume that $S=\emptyset$ so that 
$$
f= a(\{i,k\})z_iz_k+a(\{i,\ell\})z_iz_\ell + a(\{j,k\})z_jz_k+a(\{j,\ell\})z_jz_\ell + a(\{\ell,k\})z_\ell z_k
$$
Then 
$$
\Delta_{ij}(f)=  \Big(a(\{i,k\})z_k+a(\{i,\ell\})z_\ell\Big)\Big(a(\{j,k\})z_k+a(\{j,\ell\})z_\ell\Big).  
$$
It follows that $\Delta_{ij}(f)(x_\ell,x_k) \geq 0$ for all 
$x_\ell,x_k \in \RR$ if and only if 
$$
a(\{i,k\})a(\{j,\ell\})=a(\{i,\ell\})a(\{j,k\}),
$$
which proves the lemma by Theorem \ref{logf}. 
\end{proof}
\begin{lemma}\label{konstig}
Suppose that $\{i,j,x\},\{i,j,y\},\{i,k,x\},\{i,k,y\}$ are different bases of $F_7$. Then $\{i,j,k\} \notin F_7$ or 
$\{i,y,x\} \notin F_7$. 
\end{lemma}
\begin{proof}
If $\{i,j,k\} \notin F_7$ then the conclusion holds, so assume instead that $\{i,j,k\} \in F_7$. If $F_7$ has lines $\{i,j,p\} \notin F_7$ and $\{i,k,q\} \notin F_7$ then no two of $i,j,k,p,q,x,y$ are equal. It follows that 
$\{i,x,y\} \notin F_7$ is a line.  

\end{proof}
Fix $1 \leq x<y \leq 7$ and let $G_{xy}$ be the graph with vertex set 
$$
V= \{ \{i,j\} : \{i,j,x\}, \{i,j,y\} \in \mathcal{B}(F_7) \}
$$ 
and edges between sets that have non-empty intersection. 
Then $G_{xy}$ is connected, see Figure~\ref{G}. 
\setlength{\unitlength}{14mm}
\begin{figure}\caption{\label{G} The graph  $G_{67}$.}
\begin{center}
\begin{picture}(8,3.5)
\put(3,0){
\put(1,3){\p}\put(1.1,3){\{4,5\}}
\put(0,2){\p}\put(-0.7,2){\{3,4\}}\put(2,2){\p}\put(2.1,2){\{3,5\}}
\put(0,1){\p}\put(-0.7,1){\{2,3\}}\put(2,1){\p}\put(2.1,1){\{1,3\}}
\put(1,0){\p}\put(1.1,-0.1){\{1,2\}}
\path(1,3)(0,2)(0,1)(1,0)(2,1)(2,2)(1,3)
\path(0,1)(2,1)(0,2)(2,2)(0,1)
}
\end{picture}
\end{center}
\end{figure}
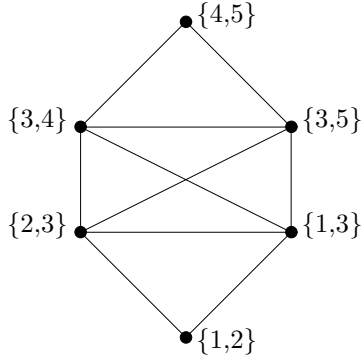
Assume now that we have a real stable polynomial $f=\sum_{S \in F_7}a(S)z^S$ with $\supp(f)=F_7$. 
The remainder of this section is devoted to deriving a contradiction under this assumption. Lemma~\ref{key} gives the following relations between the coefficients: Let $a,b,c$ be on a 
line and $d,e \in \{1,\ldots, 7\} \setminus \{a,b,c\}$. Then 
\begin{equation}\label{relations}
\frac {a(\{a,b,d\})}{a(\{a,b,e\})} = \frac {a(\{a,c,d\})}{a(\{a,c,e\})} = \frac {a(\{b,c,d\})}{a(\{b,c,e\})}.
\end{equation}

\begin{lemma}
Let $1 \leq x<y \leq 7$. 
The quotient 
$$
\frac {a(\{i,j,x\})} {a(\{i,j,y\})}
$$
is the same for all $\{i,j\}$ in $G_{xy}$.
\end{lemma}
\begin{proof}
Since $G_{xy}$ is connected it suffices to prove that 
\begin{equation}\label{kvot}
\frac {a(\{i,j,x\})} {a(\{i,j,y\})}=\frac {a(\{i,k,x\})} {a(\{i,k,y\})}
\end{equation}
whenever $\{i,j,x\},\{i,j,y\},\{i,k,x\},\{i,k,y\} \in F_7$. By Lemma \ref{konstig} we have either 
$\{i,j,k\} \notin F_7$ or $\{i,y,x\} \notin F_7$. In the first case we get by letting $i,j,k,x,y =a,b,c,d,e$ in 
\eqref{relations} that 
$$
\frac {a(\{i,j,x\})} {a(\{i,j,y\})}=\frac {a(\{i,k,x\})} {a(\{i,k,y\})}.
$$
 In the second case we have by letting $i,y,x,j,k = a,b,c,d,e$ in 
\eqref{relations} that 
$$
\frac {a(\{i,y,j\})} {a(\{i,y,k\})}=\frac {a(\{i,x,j\})} {a(\{i,x,k\})}, 
$$
which is equivalent to \eqref{kvot}.
\end{proof}
For two distinct numbers $x,y \in \{1, \ldots, 7\}$ let $\lambda_{xy} = a(\{i,j,x\})/a(\{i,j,y\})$, where $\{i,j\} \in G_{xy}$. 

\begin{lemma}\label{lambda}
Let $x,y,z \in \{1,\ldots, 7\}$ be distinct. Then 
$$
\lambda_{xz}=\lambda_{xy}\lambda_{yz}
$$
\end{lemma}
\begin{proof}
If $\{x,y,z\} \in F_7$ then there are $i,j$ such that $\{i,j,x\}, \{i,j,y\},\{i,j,z\} \in F_7$. Hence, 
$$
\lambda_{xz}= \frac{a(\{i,j,x\})}{a(\{i,j,y\})} \cdot \frac{a(\{i,j,y\})}{a(\{i,j,z\})} = \lambda_{xy}\lambda_{yz}.
$$
If $\{x,y,z\} \notin F_7$ then  for all $u \notin \{x,y,z\} $ we have $\{x,u,z\}, \{u,y,z\} \in F_7$. Hence 
\begin{eqnarray*}
\lambda_{xz} &=& \lambda_{xu}\lambda_{uz} \\
                          &=& \lambda_{xy}\lambda_{yu}\lambda_{uy}\lambda_{yz}\\
                          &=& \lambda_{xy}\lambda_{yz}.
\end{eqnarray*}
\end{proof} 

\begin{lemma}\label{po}
There are positive numbers $v_i$, $1\leq i \leq 7$ and a complex number $C$ such that 
$$
a(\{i,j,k\})= Cv_iv_jv_k
$$
for all $\{i,j,k\} \in F_7$. 
\end{lemma}
\begin{proof}
Let $v_i=\lambda_{i1}$, so that $\lambda_{xy}= v_x /v_y$ for all $1\leq x,y \leq 7$. Let $A=\{i,j,k\} \in F_7,B= \{\ell,m,n\} \in F_7$. If $i=\ell$ and $j=m$ then 
$$
\frac {a(A)}{a(B)} = \lambda_{kn}= \frac{v_iv_jv_k}{v_\ell v_m v_n}. 
$$
Otherwise, by the exchange axiom, there is 
a path 
$$
A=A_1 \rightarrow A_2 \rightarrow \cdots \rightarrow A_p = B 
$$
such that $|A_i\cap A_{i+1}|=2$. Hence,  
$$
\frac {a(A)}{a(B)}= \frac {a(A_1)}{a(A_2)} \cdots \frac {a(A_{p-1})}{a(A_{p})} =  
\frac{v_iv_jv_k}{v_\ell v_m v_n}. 
$$
Consequently 
$$
\frac {a(i,j,k)}{v_iv_jv_k} =C
$$
does not depend on $i,j,k$. 
\end{proof}
\begin{theorem}
There is no stable polynomial whose support is $F_7$.
\end{theorem}
\begin{proof}
If there were such a polynomial then by the change of variables $z_i \mapsto z_i/v_i$ and 
Lemma \ref{po} we would have that 
$$
\sum_{S \in F_7} z^S
$$
is stable. This is not the case, see \cite{COSW}.
\end{proof}

%



\section{Open Problems}
Can  the technique in Section \ref{fanosec} be extended to other matroids besides the Fano matroid? In particular,  can this technique be used to prove that the non-Pappus matroid is not the support of a polynomial with the half-plane property?  Even better, can we characterize the matroids and jump systems that are supports of polynomials with the half-plane property in matroid theory terms? 
\begin{question}
Suppose that  $f \in \CC[z_1, \ldots, z_n]$ is stable. Is there a Hermitian matrix $H$ and 
positive semidefinite matrices $A_i$ such that 
$$  
\supp(f)= \supp(\det(z_1A_1 + \cdots + z_nA_n + H))?
$$
\end{question}
It is not true that every real stable polynomial in $n \geq 3$ variables can be written as 
$\det(z_1A_1 + \cdots + z_nA_n + H)$. However, it is likely that the class of such polynomials is large enough so that all supports of stable polynomials are supports of determinants of pencils of matrices.  
\begin{question}
Are all finite jump systems in $\NN^2$ supports of stable polynomials?
\end{question}
When looking for $H$-stable polynomials with a given support it is enough to look among the 
real stable polynomials. 
\begin{proposition}
Let $f$ be a polynomial with the half-plane property. Then there is a real stable polynomial $\tilde{f}$ 
with 
$$
\supp(f)=\supp(\tilde{f}).
$$
\end{proposition}
\begin{proof}
By a rotation of the variables we may assume that $f=h+ig$, $h,g \in \RR[z_1, \ldots, z_n]$ is stable. By 
Corollary \ref{addz} we $h+z_{n+1}g$ is real stable, so $h+\alpha g$ is real stable for 
every $\alpha \in \RR$. Also, 
$$
\supp(h+\alpha g) = \supp(f), 
$$
for all but finitely many $\alpha \in \RR$. 
\end{proof}

\bigskip

\begin{center}
{  A{\Small CKNOWLEDGEMENTS}} 
\end{center}

\medskip

I would like to thank YoungBin Choe, Alan Sokal and two anonymous referees for helpful comments and suggestions, and  one of the referees for providing a cleaner proof of Lemma~\ref{konstig}.  I would also like to thank John Stembridge for helpful discussions. 
  
\end{document}